\RequirePackage{etoolbox}
\csdef{input@path}{}
\csgdef{bibdir}{}

\documentclass[ba]{imsart}
\pubyear{0000}
\volume{00}
\issue{0}
\doi{0000}
\firstpage{1}
\lastpage{1}

\usepackage{bm}
\usepackage{amsthm}
\usepackage{amsmath}
\usepackage[colorlinks,citecolor=blue,urlcolor=blue,filecolor=blue,backref=page]{hyperref}
\usepackage{natbib}
\usepackage{xspace}

\usepackage[english]{babel}
\usepackage[T1]{fontenc}
\usepackage[utf8]{inputenc}
\usepackage[scaled=1.04]{biolinum}

\usepackage{fourier}
\usepackage[scaled=0.83]{beramono}

\newcommand{\pvalue}{\textit{p-}value\xspace}
\newcommand{\given}{\, \middle| \,}
\newcommand{\pg}[2]{p\mathopen{}\left(#1 \given #2 \right)\mathclose{}}
\newcommand{\Pg}[2]{P\mathopen{}\left(#1 \given #2 \right)\mathclose{}}
\newcommand{\givenfixed}{\, | \,}
\newcommand{\pgfixed}[2]{p(#1 \givenfixed #2)}
\newcommand{\intd}{\text{d}}

\newcommand{\data}{\mathcal{D}}
\newcommand{\reject}{\mathcal{R}}
\newcommand{\param}{\theta}
\newcommand{\params}{\boldsymbol{\param}}
\newcommand{\dparam}{\intd\param}
\newcommand{\dparams}[1][n]{\intd^{#1}\param}
\newcommand{\region}[2]{\left\{#1: \, #2\right\}}

\newcommand{\refapp}{App.~\ref}
\newcommand{\refeq}{eq.~\ref}
\newcommand{\refsec}{sec.~\ref}

\usepackage[usenames]{xcolor}

\startlocaldefs
\numberwithin{equation}{section}
\theoremstyle{plain}

\endlocaldefs

\begin{document}

\begin{frontmatter}

\title{Neyman-Pearson lemma for Bayes factors}
\runtitle{Neyman-Pearson lemma for Bayes factors}

\begin{aug}

\author{\fnms{Andrew} \snm{Fowlie}\thanksref{addr1}\ead[label=e1]{andrew.j.fowlie@googlemail.com}}

\runauthor{A.~Fowlie}

\address[addr1]{Department of Physics and Institute of Theoretical Physics, Nanjing Normal University,
  Nanjing, Jiangsu 210023, China
  \href{mailto:andrew.j.fowlie@njnu.edu.cn}{andrew.j.fowlie@njnu.edu.cn} 
}

\end{aug}

\begin{abstract}
We point out that the Neyman-Pearson lemma applies to Bayes factors if we consider \emph{expected} type-1 and type-2 error rates.
That is, the Bayes factor is the test statistic that maximises the expected power for a fixed expected type-1 error rate.
For Bayes factors involving a simple null hypothesis, the expected type-1 error rate is just the completely frequentist type-1 error rate.
%
%
Lastly we remark on connections between the Karlin-Rubin theorem and uniformly most powerful tests, and Bayes factors.
This provides frequentist motivations for computing the Bayes factor and could help reconcile Bayesians and frequentists.
\end{abstract}
\end{frontmatter}

\section{Introduction}

Testing models and hypotheses against experimental data is a fundamental part of science yet the statistical approaches for doing it remain contentious. In many fields null hypothesis significance testing via \pvalue{}s is popular. In this approach, one tests a particular null hypothesis by computing the probability of obtaining data as or more extreme than that observed, assuming that the null hypothesis were true. This probability, the \pvalue, may be used as either a measure of evidence against the null hypothesis, or as a means to controlling a type-1 error rate~\citep{10.2307/20115217}. This error rate is the rate at which we would reject the null hypothesis when it was true. If we reject the null hypothesis only when $p < \alpha$, the type-1 error rate would be $\alpha$.

The use of \pvalue{}s has been contentious for decades~\citep{Edwards1992}, and in the recent years the situation appeared to reach a critical point~\citep{benjamin2018redefine,lakens2018justify}, after it became apparent that many effects discovered through null hypothesis significance testing were in fact spurious and could not be reproduced by other researchers~\citep{10.1371/journal.pmed.0020124}. All kinds of remedies have been proposed, from abandoning null hypothesis significance testing altogether~\citep{mcshane2019abandon} through to Bayesian methods, including Bayes factors~\citep{Jeffreys:1939xee,doi:10.1080/01621459.1995.10476572}.

In the Bayesian approach we compute directly the relative probability that two models, $H_0$ and $H_1$, predict the observed data, $\data$,
\begin{equation}\label{eq:bf}
B(\data) \equiv \frac{\pg{\data}{H_1}}{\pg{\data}{H_0}}.
\end{equation}
This so-called Bayes factor updates the relative plausibility of two models. The numerator and denominator in the Bayes factor are Bayesian evidences. In general for composite models they may be written
\begin{equation}
\pg{\data}{H} = \int \pg{\data}{H, \params} \pg{\params}{H}\,\dparams,
\end{equation}
that is, as the marginal likelihood, where $\params$ are the model's $n$ unknown parameters with prior density $\pg{\params}{H}$.

In the null hypothesis significance test, we are left to decide what possible data would be more extreme than that observed. If we consider a fixed type-1 error rate, this amounts to choosing a rejection region, $\reject$, in the sample space of the experiment of size $\alpha$. Although other points of view exist, the dominant approach is choosing a rejection region that maximises the chance of rejecting the null hypothesis when it is false, i.e., maximising the power, $1 - \beta$. In other words, we chose $\reject$ such that for fixed
\begin{equation}\label{eq:alpha}
\alpha \equiv \Pg{\reject}{H_0}
\end{equation}
we maximise
\begin{equation}
1 - \beta \equiv \Pg{\reject}{H_1}.
\end{equation}
Equivalently, we minimise the type-2 error rate $\beta$, the probability of not rejecting the null hypothesis when it is false. The Neyman-Pearson lemma~\citep{doi:10.1098/rsta.1933.0009} and the Karlin-Rubin theorem~\citep{casella2002statistical} indicate the most powerful tests are often built upon likelihood ratios, and likelihood ratios are ubiquitous in significance tests.

We could, on the other hand, relinquish control of the type-1 error rate, but insist on a rejection region that minimises a weighted sum of type-1 and type-2 error rates~\citep{Lindley1953,Lehmann1958,MR0133898,10.1371/journal.pone.0032734},
\begin{equation}\label{eq:weighted_sum}
W \equiv w_1 \alpha + w_2 \beta,
\end{equation}
with weights $w_1$ and $w_2$. For the purposes of the calculations and arguments in this work, this is equivalent to the conventional approach discussed above. To see that, suppose that the region that minimises \refeq{eq:weighted_sum} has size $\alpha^\prime$. Since this region minimises $W$, among all possible regions of size $\alpha^\prime$, it must be the one that maximises the power, $1 - \beta$. 

In this work, we use the Bayes factor as a test statistic in a frequentist setting, argue that it is an optimal choice and that this may help reconcile advocates of Bayesian and frequentist procedures for hypothesis testing. The idea of using the Bayes factor in this way and its potential for harmonising Bayesian and frequentist procedures~\citep{berger2003,BAYARRI201690} dates as far back as the 1950s~\citep{10.1214/aoms/1177706790,good1961weight,10.2307/2290192}, though \cite{10.1214/ss/1030037904} notes that it failed to gain traction among practitioners. In \refsec{sec:np} we briefly review the Neyman-Pearson lemma and the Karlin-Rubin theorem, and in \refsec{sec:np_for_bf} and \refsec{sec:kr_for_bf}, we show they may be applied to Bayes factors, such that the Bayes factor is the most powerful choice of test statistic. This was previously discussed formally in \cite{doi:10.1111/anzs.12171}. Besides drawing attention to the relevant results and presenting them in an original, pedagogical manner, we combine the ideas of~\cite{10.1214/aoms/1177706790,good1961weight,10.2307/2290192} and results of~\cite{doi:10.1111/anzs.12171}. We thus argue in \refsec{sec:conclusions} that the Neyman-Pearson lemma and Karlin-Rubin theorem in fact put Bayes factors at the centre of frequentist hypothesis tests as they are an optimal choice of test statistic. This may help finally reconcile Bayesian and frequentist approaches to testing.

\section{The Neyman-Pearson lemma}\label{sec:np}

For simple hypotheses with no unknown parameters, the Neyman-Pearson lemma~\citep{doi:10.1098/rsta.1933.0009} tell us that the likelihood ratio,
\begin{equation}
\lambda(\data) = \frac{\pg{\data}{H_1}}{\pg{\data}{H_0}},
\end{equation}
is an optimal test statistic. The rejection region of size $\alpha$ that maximises the power must be a contour of the likelihood ratio, $\reject = \region{\data}{\lambda(\data) \ge \lambda_0}$, where the value at the contour, $\lambda_0$, depends on the desired type-1 error rate and would be found from \refeq{eq:alpha}.
The Neyman-Pearson lemma doesn't extend to composite hypotheses that depend on unknown parameters, as the optimal choice of test statistic would in general depend on the assumed values for the unknown parameters. In some cases, a uniformly most powerful test (UMPT) exists, which maximises the power for any values of the unknown parameters in the alternative hypothesis. The Karlin-Rubin theorem demonstrates cases in which a UMPT exists; see \refsec{sec:kr_for_bf}.


\section{The Neyman-Pearson lemma for Bayes factors}\label{sec:np_for_bf}

In the Bayesian formalism, for any composite model, we may find a simple model by marginalising the $n$ unknown parameters,
\begin{equation}
\pg{\data}{H} = \int \pg{\data}{H, \params} \pg{\params}{H}\,\dparams.
\end{equation}
This is the prior predictive for the data and when evaluated at the observed data, the Bayesian evidence for the model. We may thus use our Bayes factor in place of the likelihood ratio in the Neyman-Pearson lemma. In this case, though, the error rates are subtly re-interpreted, and to distinguish them we denote them with a bar,
\begin{align}\label{eq:bar_alpha}
\bar\alpha &\equiv \Pg{\reject}{H_0} = \int \Pg{\reject}{H_0, \params_0} \, \pg{\params_0}{H_0}\,\dparams[n]_0\\
1 - \bar\beta &\equiv \Pg{\reject}{H_1} = \int \Pg{\reject}{H_1, \params_1} \, \pg{\params_1}{H_1}\,\dparams[m]_1
\end{align}
That is, here they are the \emph{expected} error rates, averaging over the possible values for the $n$ unknown model parameters $\params_0$ in the null hypothesis and for the $m$ unknown model parameters $\params_1$ in the alternative hypothesis. They do not in general correspond to any observable long-run error rates.

This requires, of course, choices of prior density, $\pg{\params_0}{H_0}$ and $\pg{\params_1}{H_1}$; see e.g., \cite{doi:10.1080/01621459.1996.10477003,consonni2018} for discussion of rules for choosing priors in a Bayesian setting. In a frequentist setting, the priors needn't be interpreted in the same manner as in subjective or objective Bayesian approaches. The prior for any parameters in the alternative hypothesis, $\pg{\params_1}{H_1}$, may be thought of as a weight function that indicates which choices of parameters we want our test to be most powerful for~\citep{10.2307/2286006,BAYARRI201690}. Similarly, the prior for any parameters in a composite null hypothesis, $\pg{\params_0}{H_0}$, weights which choices of parameters we most want to control the type-1 error rate for.

We could instead consider a fixed \emph{maximum} type-1 error rate, $\hat \alpha$, for any values of the unknown parameters in the null hypothesis. If we assume that it occurs when $\params_0 = \hat\params_0$, we may write
\begin{equation}\label{eq:hat_alpha}
\hat \alpha = \Pg{\reject}{H_0, \hat\params_0}.
\end{equation}
This is in fact equivalent to a sharp prior $\pg{\params_0}{H_0} = \delta(\params_0 - \hat\params_0)$, i.e., specifying through our prior that we must control the type-1 error rate for the parameters that maximise the type-1 error rate. In this case we would use the Bayes factor,
\begin{equation}\label{eq:hat_bf}
B(\data) = \frac{
  \int \pg{\data}{H_1, \params_1} \, \pg{\params_1}{H_1} \, \dparams_1
  }{
  \pg{\data}{H_0, \hat\params_0}
  }
\end{equation}
in place of the likelihood ratio in the Neyman-Pearson lemma. Under the null hypothesis the expected magnitude of the preference for the alternative model from the Bayes factor in \refeq{eq:bf} must be smaller than that from the one used when we control the maximum type-1 error rate in \refeq{eq:hat_bf}. This result follows from Gibbs' inequality,
\begin{equation}
\int \left[\log\left(\frac{
  \pg{\data}{H_1}
}{
  \pg{\data}{H_0 \vphantom{, \hat\params_0}}
}\right)
-
\log\left(\frac{
  \pg{\data}{H_1}
}{
  \pg{\data}{H_0, \hat\params_0}
}\right)\right]
\, \pg{\data}{H_0} \, \intd \data
\le 0.
\end{equation}
So whilst controlling the maximum error rate rather than the expected error may seem conservative, the computation involves the Bayes factor in \refeq{eq:hat_bf} that we expect to overstate the evidence against the null hypothesis.

If the null hypothesis is simple, as is often the case, $\bar \alpha = \hat\alpha = \alpha$, and so the type-1 error rates may be interpreted in the usual, completely frequentist manner. The power, on the other hand, remains the expected power, averaged across the unknown parameters in the alternative hypothesis.

\section{Karlin-Rubin theorem for Bayes factors}\label{sec:kr_for_bf}

It would be desirable to extend the Neyman-Pearson lemma to composite models. Unfortunately, UMPT do not always exist as the optimal test generally depends on the values of the unknown parameters in the alternative hypothesis. There are approaches that sidestep the issue such as a minimax treatment of type-1 error rates and power. The Karlin-Rubin theorem~\citep{casella2002statistical}, on the other hand, extends the Neyman-Pearson lemma to a UMPT in a composite case in special circumstances.

To apply the theorem, the null and alternative hypothesis must be disjoint regions of a one-dimensional parameter space separated by a boundary at $\theta_C$, that is,
\begin{equation}\label{eq:kr_H}
H_0: \, \param \le \param_C \quad\text{and}\quad H_1: \, \param > \param_C.
\end{equation}
Suppose that a sufficient test statistic, $T$, exists and that the likelihood ratio
\begin{equation}\label{eq:mlr}
\lambda(T; \param, \param^\prime) = \frac{\pg{T}{\param}}{\pg{T}{\param^\prime}}
\end{equation}
is a monotonic non-decreasing function of $T$ for any $\param > \param^\prime$. We call these condition the monotone likelihood ratio (MLR) conditions. Under these conditions, the UMPT is a threshold on $T$, $\reject = \region{T}{T \ge T_0}$. The threshold $T_0$ is determined by fixing the maximum type-1 error rate,
\begin{equation}
\hat\alpha = \Pg{\reject}{H_0, \param_C}
\end{equation}
and occurs when $\param = \param_C$. This ensures the size of the test is no larger than $\hat\alpha$ for any choice of unknown parameter in the null hypothesis.

We can make a similar statement for our Bayes factor,
\begin{equation}\label{eq:kr_bf}
B(T) = \frac{\pg{T}{H_1}}{\pg{T}{H_0}}
  = \frac{
  \int \pg{T}{H_1, \params_1} \, \pg{\params_1}{H_1} \, \dparams[m]_1
  }{
  \int \pg{T}{H_0, \params_0} \, \pg{\params_0}{H_0} \, \dparams_0.
  }
\end{equation}
We assume that the prior $\pg{\params_0}{H_0}$ for the unknown parameters in the null hypothesis was chosen. We suppose that, with that choice, the Bayes factor is always a monotonic function of $T$ for any choice of prior for the $m$ unknown parameters in the alternative model, $\pg{\params_1}{H_1}$. This is satisfied by the MLR conditions of the Karlin-Rubin theorem for any choice of $\pg{\params_0}{H_0}$; see \refapp{app:kr}.
The Karlin-Rubin theorem considered a fixed maximum type-1 error rate, $\hat \alpha$. In our Bayesian interpretation, we could fix this or the mean error rate $\bar\alpha$. If we fix the former, we consider the Bayes factor in \refeq{eq:hat_bf}
in the following argument. If we choose to fix the expected type-1 error rate, on the other hand, we instead consider the Bayes factor in \refeq{eq:kr_bf} in the following argument.

By an application of the Neyman-Pearson lemma for Bayes factors, the most powerful test at a fixed size should be a threshold on $B$, $\reject = \region{T}{B(T) \ge B_0}$. As the Bayes factor is a monotonic function of $T$, this is equivalent to a threshold on $T$, $\reject = \region{T}{T \ge T_0}$. As $T_0$ can be found independently from the prior $\pg{\params_1}{H_1}$ by eq.~\ref{eq:bar_alpha} or  \ref{eq:hat_alpha}, it must be the most powerful test for any choice of prior $\pg{\params_1}{H_1}$, including point masses at any particular values of the unknown parameters. Thus, it is the UMPT.

We thus find connections between the Bayes factor and the UMPT. By generalising the Karlin-Rubin theorem, we find that a UMPT exists whenever the Bayes factor is a monotonic function of a sufficient statistic for any choice of prior for the unknown parameters in the alternative model. The Bayes factor corresponding the observed $T$ must, however, depend on the choices of prior.

\section{Reflections}\label{sec:conclusions}

The Neyman-Pearson lemma leads to optimal test statistics for null hypothesis significance tests for simple hypotheses. Interpreting the Bayes factor as a likelihood ratio for two simple models leads to a Bayesian interpretation of the Neyman-Pearson lemma that goes beyond simple models. The Bayes factor maximises the expected power for a test of a fixed expected size.

On the Bayesian side, this could provide further justification for using Bayes factors for objective Bayesians or more generally those who are concerned about the frequentist properties of Bayes factors. If we place a threshold on the Bayes factor, for whatever type-1 error rate it to which that threshold corresponds, the Bayes factor was the statistic that maximised the expected power. On the frequentist side, even if you want to carry on computing \pvalue{}s, there is justification for doing so using the Bayes factor as a test statistic, especially in the case of simple null hypotheses but composite alternatives. In the case of simple null hypotheses, using the Bayes factor results in the best \emph{expected} power for a fixed completely frequentist type-1 error rate. The only concession required is that to construct the test we must choose a weight function that marks where we want power and talk about expected power. The test itself would, however, remain strictly frequentist.

The Karlin-Rubin theorem extends the Neyman-Pearson lemma to particular composite models. We found that the conditions of the Karlin-Rubin theorem may be recast as the requirement that the Bayes factor is a monotonic function of a sufficient statistic for any choices of prior for unknown parameters. This lead to a slightly novel proof of a generalised Karlin-Rubin theorem and a connection between the properties of Bayes factors and the existence of uniformly most powerful tests. The Bayesian interpreted Karlin-Rubin theorem provides conditions under which a test of fixed size always maximises the expected power, regardless which prior was chosen for the unknown parameters in the alternative hypothesis in the computation of the expected power.

These results could help synthesise frequentist and Bayesian procedures, as it shows that the Bayes factors could lie at the heart of each one, and proponents of either should be interested, in principle, in computing the Bayes factor. The outstanding difference would be that the Bayesian would consider the magnitude of the observed Bayes factor, in accordance with the likelihood principle~\citep{zbMATH02166302}, whereas the frequentist would consider the probability of obtaining a Bayes factor more extreme than that observed. In practice, computing the Bayes factor and finding its distribution could be challenging, as popular asymptotic approaches such as Wilks' theorem~\citep{wilks1938} needn't apply (though see~\cite{10.2307/25734099}).

\appendix

\section{Karlin-Rubin conditions for Bayes factors}\label{app:kr}

Let us check the implications of the MLR conditions in the Karlin-Rubin theorem on the Bayes factor. For the hypotheses under consideration in the Karlin-Rubin theorem in \refeq{eq:kr_H}, the Bayes factor could be written
\begin{equation}
B(T) = \frac{
  \int_{\param > \param_C} \pg{T}{\param} \, \pg{\param}{H_1} \, \dparam
  }{
  \int_{\param \le \param_C} \pg{T}{\param} \, \pg{\param}{H_0} \, \dparam
  }
\end{equation}
for some suitably normalised choices of prior densities $\pg{\param}{H_1}$ and $\pg{\param}{H_0}$.
Starting from \refeq{eq:mlr}, for $T^\prime > T$ and $\param > \param^\prime$, we have
\begin{align}
\frac{\pgfixed{T^\prime}{\param}}{\pgfixed{T^\prime}{\param^\prime}} &\ge \frac{\pgfixed{T}{\param}}{\pgfixed{T}{\param^\prime}}\\
\frac{\pgfixed{T^\prime}{\param}}{\pgfixed{T^\prime}{\param^\prime}} \, \frac{\pgfixed{\param}{H_1}}{\pgfixed{\param^\prime}{H_0}}
&\ge
\frac{\pgfixed{T}{\param}}{\pgfixed{T}{\param^\prime}} \, \frac{\pgfixed{\param}{H_1}}{\pgfixed{\param^\prime}{H_0}}\\
\left[\pg{T^\prime}{\param} \, \pg{\param}{H_1}\right] \left[\pg{T}{\param^\prime} \, \pg{\param^\prime}{H_0}\right]
&\ge
\left[\pg{T}{\param} \, \pg{\param}{H_1}\right] \left[\pg{T^\prime}{\param^\prime} \, \pg{\param^\prime}{H_0}\right]
\end{align}
where we multiplied each side by the ratio of prior densities that appear in the Bayes factor and rearranged terms. Integrating each side with respect to $\param$ and $\param^\prime$ only over the regions guaranteeing $\param > \param^\prime$, we find
\begin{equation}
\begin{split}
\int_{\param > \param_C} \pg{T^\prime}{\param} \, \pg{\param}{H_1} \, \dparam \,
& \int_{\param \le \param_C} \pg{T}{\param} \, \pg{\param}{H_0} \, \dparam
\ge\\
& \int_{\param > \param_C}\pg{T}{\param} \, \pg{\param}{H_1} \, \dparam \,
\int_{\param \le \param_C} \pg{T^\prime}{\param} \, \pg{\param}{H_0} \, \dparam.
\end{split}
\end{equation}
Finally, we rearrange terms to find
\begin{equation}
\begin{split}
\frac{
  \int_{\param > \param_C} \pg{T^\prime}{\param} \, \pg{\param}{H_1} \, \dparam
}{
  \int_{\param \le \param_C} \pg{T^\prime}{\param} \, \pg{\param}{H_0} \, \dparam
}
\ge
\frac{
  \int_{\param > \param_C} \pg{T}{\param} \, \pg{\param}{H_1} \, \dparam
}{
  \int_{\param \le \param_C} \pg{T}{\param} \, \pg{\param}{H_0} \, \dparam
}
\end{split}
\end{equation}
such that $B(T^\prime) \ge B(T)$ for $T^\prime \ge T$. We did not make use of any properties of any particular choice of prior densities. In other words, the ordinary MLR conditions of the Karlin-Rubin theorem mean that the Bayes factor is a monotonic function of $T$ for any choices of prior densities $\pg{\param}{H_1}$ and $\pg{\param}{H_0}$.

\bibliographystyle{ba}
\bibliography{references}
\end{document}